\documentclass[12pt]{article}
\usepackage{latexsym}
\usepackage{amssymb}
\usepackage{graphicx}
\usepackage{cite}

\newtheorem{Theorem}{Theorem}[part]
\newtheorem{Definition}{Definition}[part]
\newtheorem{Proposition}{Proposition}[part]

\newtheorem{Lemma}{Lemma}[part]
\newtheorem{Corollary}{Corollary}[part]
\newtheorem{Remark}{Remark}[part]
\newtheorem{Example}{Example}[part]
\newtheorem{Algorithm}{Algorithm}

\def \ep{\hbox{ }\hfill$\Box$}

\def \ra{\rightarrow}
\def\reff#1{{\rm(\ref{#1})}}

\begin{document}
\title{Finding the Spectral Radius of a Nonnegative Tensor}

\author{
Shenglong Hu \thanks{Email: shenglong@tju.edu.cn. Department of
Applied Mathematics, The Hong Kong Polytechnic University, Hung Hom,
Kowloon, Hong Kong.},\hspace{4mm} Zheng-Hai Huang \thanks{Email:
huangzhenghai@tju.edu.cn. Department of Mathematics, School of
Science, Tianjin University, Tianjin, China. This author's work was
supported by the National Natural Science Foundation of China (Grant No. 10871144
and Grant No. 30870713).},\hspace{4mm} Liqun Qi \thanks{Email:
maqilq@polyu.edu.hk. Department of Applied Mathematics, The Hong
Kong Polytechnic University, Hung Hom, Kowloon, Hong Kong. This
author's work was supported by the Hong Kong Research Grant
Council.} }

\date{\today} \maketitle

\begin{abstract}
In this paper, we introduce a new class of nonnegative tensors --- strictly nonnegative tensors.
A weakly irreducible
nonnegative tensor is a strictly nonnegative tensor but not vice
versa.   We show that the spectral radius of a strictly nonnegative
tensor is always positive. We give some sufficient and necessary
conditions for the six well-conditional classes of nonnegative
tensors, introduced in the literature, and a full relationship
picture about strictly nonnegative tensors with these six classes of
nonnegative tensors. We then establish global R-linear convergence
of a power method for finding the spectral radius of a nonnegative
tensor under the condition of weak irreducibility. We show that for
a nonnegative tensor $T$, there always exists a partition of the
index set such that every tensor induced by the partition is weakly
irreducible; and the spectral radius of $T$ can be obtained from
those spectral radii of the induced tensors. In this way, we develop
a convergent algorithm for finding the spectral radius of {\it a
general nonnegative tensor} without any additional assumption. The
preliminary numerical results demonstrate the feasibility and
effectiveness of the proposed algorithm. \vspace{3mm}

\noindent {\bf Key words:}\hspace{2mm} Nonnegative tensor, spectral
radius, strict nonnegativity, weak irreducibility, algorithm
\vspace{3mm}

\noindent {\bf AMS subjection classifications (2010):}\hspace{2mm}
15-02; 15A18; 15A69; 65F15
\end{abstract}

\newpage
\section{Introduction}
\setcounter{Assumption}{0}
\setcounter{Theorem}{0} \setcounter{Proposition}{0}
\setcounter{Corollary}{0} \setcounter{Lemma}{0}
\setcounter{Definition}{0} \setcounter{Remark}{0}
\setcounter{Algorithm}{0}  \setcounter{Example}{0}

\hspace{4mm}
Recently, the research topic on eigenvalues of nonnegative tensors
attracted much attention \cite{cpz, cpz2, fgh, lzi, nqz, p, yy, zq,
zqx}.   Researchers studied the Perron-Frobenius theorem for
nonnegative tensors and algorithms for finding the largest
eigenvalue, i.e., the spectral radius,  of a nonnegative tensor. Six
well-conditional classes of nonnegative tensors have been
introduced: irreducible nonnegative tensors \cite{cpz}, essentially
positive tensors \cite{p}, primitive tensors \cite{cpz2}, weakly
positive tensors \cite{zqx}, weakly irreducible nonnegative tensors
\cite{fgh} and weakly primitive tensors \cite{fgh}.   Zhang, Qi and
Xu \cite{zqx} concluded the relationships among the first four
classes of nonnegative tensors. Friedland, Gaubert and Han
\cite{fgh} introduced weakly irreducible nonnegative tensors and
weakly primitive tensors.  These two classes, as their names
suggest, are broader than the classes of irreducible nonnegative
tensors and primitive tensors respectively.

In the next section, we propose a new class of nonnegative tensors, we call them {\em strictly nonnegative} tensors.
We show that the class of strictly nonnegative tensors strictly contains the class of weakly irreducible nonnegative tensors mentioned above.
We also prove that the spectral radius of a strictly nonnegative tensor is always positive. This further strengthens the Perron-Frobenius results for nonnegative tensors in the literature \cite{cpz,cpz2,fgh,yy}.

In Section 3, we give sufficient and necessary conditions for the
six well-conditional classes of nonnegative tensors, introduced in
the literature, and a full relationship picture about strictly
nonnegative tensors with these six classes of nonnegative tensors.

Friedland, Gaubert and Han \cite{fgh} proposed a power method for
finding the largest eigenvalue of a weakly irreducible nonnegative
tensor, and established its R-linear convergence under the condition
of weak primitivity. In Section 4, we modify that method and
establish its global R-linear convergence for weakly irreducible
nonnegative tensors.

Then, in Section 5, we show that for a nonnegative tensor $T$,
always there exists a partition of the index set $\{1,\ldots,n\}$
such that every tensor induced by the partition is weakly
irreducible; and the largest eigenvalue of $T$ can be obtained from
those largest eigenvalues of the induced tensors.  In Section 6,
based on the power method for weakly irreducible nonnegative tensors
proposed in Section 4, we develop a convergent algorithm for finding
the spectral radius of {\it a general nonnegative tensor} without
any additional assumption.  We report some preliminary numerical
results of the proposed method for general nonnegative tensors.
These numerical results demonstrate the feasibility and
effectiveness of the proposed algorithm. Conclusions and remarks are
given in Section 7.

\medskip

Here is some notation in this paper. A tensor $T$ in real field
$\Re$ of order $m$ and dimension $n$ with $m,n\geq 2$ is an $m$-way
array which can be denoted by $(T_{i_1\ldots i_m})$ with
$T_{i_1\ldots i_m}\in\Re$ for all $i_j\in\{1,\ldots,n\}$ and
$j\in\{1,\ldots,m\}$. For a tensor $T$ of order $m\geq 2$ and
dimension $n\geq 2$, if there exist $\lambda\in {\cal C}$ and $x\in
{\cal C}^n\setminus\{0\}$ such that
\begin{eqnarray}\label{eig}
(Tx^{m-1})_i:=\sum_{i_2,\ldots,i_m=1}^nT_{ii_2\ldots
i_m}x_{i_2}\cdots x_{i_m}=\lambda x_i^{m-1},\;\;\forall
i\in\{1,\ldots,n\}
\end{eqnarray}
holds, then $\lambda$ is called an eigenvalue of $T$, $x$ is called
a corresponding eigenvector of $T$ with respect to $\lambda$, and
$(\lambda, x)$ is called an eigenpair of $T$. This definition was
introduced by Qi \cite{q} when $m$ is even and $T$ is symmetric
(i.e., $T_{j_1\ldots j_m}=T_{i_1\ldots i_m}$ among all the
permutations $(j_1,\ldots,j_m)$ of $(i_1,\ldots,i_m)$).
Independently, Lim \cite{l} gave such a definition but restricted
$x$ to be a real vector and $\lambda$ to be a real number.  Let
$\Re^n_+:=\{x\in\Re^n \;|\;x\geq0\}$ and $\Re^n_{++}:=\{x\in\Re^n
\;|\;x>0\}$.   Suppose that $T$ is a nonnegative tensor, i.e., that
all of its entries are nonnegative. It can be seen that if we
define function $F_T:\Re^n_+\ra\Re^n_+$ associated nonnegative
tensor $T$ as
\begin{eqnarray}\label{map-1}
(F_T)_i(x):=\left(\sum_{i_2,\ldots, i_m=1}^nT_{ii_2\ldots
i_m}x_{i_2}\cdots x_{i_m}\right)^{\frac{1}{m-1}}
\end{eqnarray}
for all $i\in\{1,\ldots,n\}$ and $x\in\Re^n_+$, then \reff{eig} is
strongly related to the eigenvalue problem for the nonlinear map $F_T$
discussed in \cite{n}. Denote by
$\rho(T):=\max\{|\lambda|\;|\;\lambda\in \sigma(T)\}$ where
$\sigma(T)$ is the set of all eigenvalues of $T$. We call $\rho(T)$
and $\sigma(T)$ the spectral radius and spectra of the tensor $T$,
respectively. Hence, the eigenvalue problem considered in this paper
can be stated as: {\em for a general nonnegative tensor $T$, how
to find out $\rho(T)$}?

\section{Strictly nonnegative tensors}
\setcounter{Assumption}{0}
\setcounter{Theorem}{0} \setcounter{Proposition}{0}
\setcounter{Corollary}{0} \setcounter{Lemma}{0}
\setcounter{Definition}{0} \setcounter{Remark}{0}
\setcounter{Algorithm}{0}  \setcounter{Example}{0}
\hspace{4mm} In this section, we propose and analyze a new class of
nonnegative tensors, namely {\em strictly nonnegative tensors}. To
this end, we first recall several concepts related to nonnegative
tensors in the literature \cite{cpz,cpz2,fgh,p,zqx}.
\begin{Definition}\label{def-basic}
Suppose that $T$ is a nonnegative tensor of order $m$ and dimension
$n$.
\begin{itemize}
\item $T$ is called {\em reducible} if there exists a nonempty proper
index subset $I\subset \{1,\ldots,n\}$ such that
\begin{eqnarray}\label{red}
T_{i_1i_2\ldots i_m}=0,\quad \forall i_1\in I,\quad \forall
i_2,\ldots,i_m\notin I.
\end{eqnarray}
If $T$ is not reducible, then $T$ is called {\em irreducible}.
\item $T$ is called {\em essentially positive} if $Tx^{m-1}\in\Re^n_{++}$ for any nonzero
$x\in\Re^n_+$.
\item $T$ is called {\em primitive} if for some positive integer $k$, $F_T^{k}(x)\in\Re^n_{++}$ for any nonzero
$x\in\Re^n_+$, here $F_T^{k}:=F_T(F_T^{k-1})$.
\item A nonnegative matrix $M(T)$ is called the {\em majorization}
associated to nonnegative tensor $T$, if the $(i,j)$-th element of
$M(T)$ is defined to be $T_{ij\ldots j}$ for any
$i,j\in\{1,\ldots,n\}$. $T$ is called {\em weakly positive} if
$\left[M(T)\right]_{ij}>0$ for all $i\neq j$.
\end{itemize}
\end{Definition}

In Definition \ref{def-basic}, the concepts of {\it reducibility}
and {\it irreducibility} were proposed by Chang, Pearson and Zhang
\cite[Definition 2.1]{cpz} (an equivalent definition can be found in
Lim \cite[Page 131]{l}); the concept of {\it essential positivity} was given
in \cite[Definition 3.1]{p}; and the concept of {\it primitivity}
was given in \cite[Definition 2.6]{cpz2} while we used its
equivalent definition \cite[Thoerem 2.7]{cpz2}. In addition, the
concept of {\it majorization} was given in \cite[Definition
3.5]{cpz} and earned this name in \cite[Definition 2.1]{p}; and the
concept of {\it weak positivity} was given in \cite[Definition
3.1]{zqx}.

Friedland, Gaubert and Han \cite{fgh} defined {\it weakly
irreducible} polynomial maps and {\it weakly primitive} polynomial
maps by using the strong connectivity of a graph associated with a
polynomial map. Their concepts for homogeneous polynomials gave the
corresponding classes of nonnegative tensors.
\begin{Definition}\label{def-new}
Suppose that $T$ is a nonnegative tensor of order $m$ and dimension
$n$.
\begin{itemize}
\item We call a nonnegative matrix $G(T)$ the {\em representation}
associated to the nonnegative tensor $T$, if the $(i,j)$-th element of
$G(T)$ is defined to be the summation of $T_{ii_2\ldots i_m}$ with
indices $\{i_2,\ldots,i_m\}\ni j$.
\item We call the tensor $T$ {\em weakly reducible} if its representation $G(T)$ is a reducible matrix,
and {\em weakly primitive} if $G(T)$ is a primitive matrix. If $T$ is not
{\em weakly reducible}, then it is called {\em weakly irreducible}.
\end{itemize}
\end{Definition}

Now, we introduce {\em strictly nonnegative tensors}.

\begin{Definition}\label{add-def-2}
Suppose that $T$ is a nonnegative tensor of order $m$ and dimension
$n$. Then, it is called
{\em strictly nonnegative} if $F_T(x)>0$ for any $x>0$.
\end{Definition}

\begin{Lemma}\label{add-lem-1}
An $m$-th order $n$ dimensional nonnegative tensor $T$ is strictly
nonnegative if and only if the vector $R(T)$ with its $i$-th element
being $\sum_{i_2,\ldots,i_m=1}^nT_{ii_2\cdots i_m}$ is positive.
\end{Lemma}

\noindent {\bf Proof.} By Definition \ref{add-def-2}, $Te^{m-1}>0$ with $e$ being the vector of all ones. So,
\begin{eqnarray*}
\left(Te^{m-1}\right)_i=\sum_{i_2,\ldots,i_m=1}^nT_{ii_2\cdots i_m}>0
\end{eqnarray*}
for all $i\in\{1,\ldots,n\}$. The ``only if " part follows.

Now, suppose $R(T)>0$. Then, for every $i\in\{1,\ldots,n\}$, we could find $j_{2_i},\ldots,j_{m_i}$ such that $T_{ij_{2_i}\cdots j_{m_i}}>0$.
So, for any $x>0$, we have
\begin{eqnarray*}
\left(Tx^{m-1}\right)_i=\sum_{j_2,\ldots,j_m=1}^nT_{ij_2\cdots j_m}x_{j_2}\cdots x_{j_m}\geq T_{ij_{2_i}\cdots j_{m_i}}x_{j_{2_i}}\cdots x_{j_{m_i}}>0
\end{eqnarray*}
for all $i\in\{1,\ldots,n\}$. Hence, the ``if" part follows. The proof is complete. \ep

\begin{Corollary}\label{add-cor-1}
An $m$-th order $n$ dimensional nonnegative tensor $T$ is strictly
nonnegative if it is weakly irreducible.
\end{Corollary}

\noindent {\bf Proof.} Suppose that $T$ is weakly irreducible. Note that the signs of vectors $G(T)e$ and $R(T)$ are the same, and also that $G(T)e$ is positive since $T$ is weakly irreducible. Because, otherwise, we would have a zero row of matrix $G(T)$, which further implies that $G(T)$ is reducible, a contradiction. So, by Lemma \ref{add-lem-1}, $T$ is strictly nonnegative. \ep

A result similar to this corollary was given in (3.2) of \cite{fgh}.
On the other hand, the converse of Corollary \ref{add-cor-1} is not
true in general.
\begin{Example}\label{add-exm-1}
Let third order two dimensional nonnegative tensor $T$ be defined as:
\begin{eqnarray*}
T_{122}=T_{222}=1,\;\mbox{and}\;T_{ijk}=0\,\;\mbox{for other}\;i,j,k\in\{1,2\},
\end{eqnarray*}
then $R(T)=\left(\begin{array}{c}1\\1\end{array}\right)>0$. So $T$ is strictly nonnegative by Lemma \ref{add-lem-1}. While, $G(T)=\left(\begin{array}{cc}0&1\\0&1\end{array}\right)$ is a reducible nonnegative matrix, so $T$ is weakly reducible.
\end{Example}

\begin{Proposition}\label{add-prop-1}
An $m$-th order $n$ dimensional nonnegative tensor $T$ is strictly
nonnegative if and only if $F_T$ is strictly increasing \cite{cpz2},
i.e., $F_T(x)>F_T(y)$ for any $x>y\geq 0$.
\end{Proposition}

\noindent {\bf Proof.} If $T$ is strictly increasing, then, for any
$x>0$, we have $Tx^{m-1}>T0^{m-1}=0$. So, $T$ is strictly
nonnegative. The ``if" part follows.

Now, suppose that $T$ is strictly nonnegative. Then, $R(T)>0$ by Lemma \ref{add-lem-1}. So, for every $i\in\{1,\ldots,n\}$, we could find $j_{2_i},\ldots,j_{m_i}$ such that $T_{ij_{2_i}\cdots j_{m_i}}>0$. If $0\leq x<y$, then $x_{j_{2_i}}\cdots x_{j_{m_i}}<y_{j_{2_i}}\cdots y_{j_{m_i}}$ for any $i\in\{1,\ldots,n\}$. Hence,
\begin{eqnarray*}
\left(Ty^{m-1}\right)_i-\left(Tx^{m-1}\right)_i&=&\sum_{j_2,\ldots,j_m=1}^nT_{ij_2\cdots j_m}\left(y_{j_2}\cdots y_{j_m}-x_{j_2}\cdots x_{j_m}\right)\\
&\geq&T_{ij_{2_i}\cdots j_{m_i}}\left(y_{j_{2_i}}\cdots y_{j_{m_i}}-x_{j_{2_i}}\cdots x_{j_{m_i}}\right)\\
&>&0
\end{eqnarray*}
for any $i\in\{1,\ldots,n\}$. So, the ``only if" part follows. The
proof is complete. \ep

We now show that the spectral radius of a strictly nonnegative
tensor is always positive. At first, we present some notation. For
any nonnegative tensor $T$ of order $m$ and dimension $n$ and a
nonempty subset $I$ of $\{1,\ldots,n\}$, the induced tensor denoted
by $T_I$ of $I$ is defined as the $m$-th order $|I|$ dimensional
tensor $\{T_{i_1\ldots i_m}\;|\;i_1,\ldots,i_m\in I\}$. Here $|I|$
denotes the cardinality of the set $I$.

\begin{Lemma}\label{add-lem-2}
For any $m$-th order $n$ dimensional nonnegative tensor $T$ and nonempty subset $I\subseteq \{1,\ldots,n\}$, $\rho(T)\geq\rho(T_I)$
\end{Lemma}

\noindent {\bf Proof.} Let $K$ be a nonnegative tensor of the same size of $T$ with $K_{I}=T_{I}$ and zero anywhere else. Then, obviously, $T\geq K\geq 0$ in the sense of componentwise and $\rho(K)=\rho(T_{I})$. Now, by \cite[Lemma 3.4]{yy}, $\rho(K)\leq\rho(T)$. So, the result follows immediately.
\ep

\begin{Theorem}\label{add-thm-1}
If nonnegative tensor $T$ is strictly nonnegative, then $\rho(T)>0$.
\end{Theorem}

\noindent {\bf Proof.} By Lemma \ref{add-lem-1}, $R(T)>0$, so $G(T)e>0$. Now,
\begin{itemize}
\item [(I)] if $G(T)$ is an irreducible matrix (i.e., $T$ is a weakly irreducible tensor by Definition \ref{def-new}), then $\rho(T)$ is positive. Actually, we could find $x>0$ such that $Tx^{m-1}=\rho(T)x^{[m-1]}$ by Perron-Frobenius Theorem \cite{fgh}, so the strict positivity of $T$ implies $\rho(T)>0$.
\item [(II)] if $G(T)$ is a reducible matrix (i.e., $T$ is a weakly reducible tensor by Definition \ref{def-new}), we could find a nonempty $I\subseteq \{1,\ldots,n\}$ such that $\left[ G(T)\right ]_{ij}=0$ for all $i\in I$ and $j\notin I$. Denote by $K$ the principal submatrix of $G(T)$ indexed by $I$, and $T^{\prime}$ the tensor induced by $I$. Since $G(T)e>0$ and $\left[ G(T)\right ]_{ij}=0$ for all $i\in I$ and $j\notin I$, we still have $Ke>0$. Hence, $T^{\prime}$ is also strictly nonnegative since $G(T^{\prime})=K$ and $Ke>0$. By Lemma \ref{add-lem-2}, we have that $\rho(T)\geq\rho(T^{\prime})$.
\end{itemize}

So, inductively, we could finally get a tensor
sequence $T,T^{\prime},\ldots,T^*$ (since $n$ is finite) with
\begin{eqnarray*}
\rho(T)\geq\rho(T^{\prime})\geq\cdots\geq\rho(T^*),
\end{eqnarray*}
and $T^*$ is a weakly irreducible tensor when the dimension of $T^*$
is higher than $1$, or $T^*$ is a positive one dimensional tensor
(i.e., a scalar) since $T^*$ is strictly positive. In both cases,
$\rho(T^*)>0$ by (I). The proof is complete. \ep

\begin{Example}\label{add-exm-2}
Let third order $2$ dimensional nonnegative tensor $T$ be defined as:
\begin{eqnarray*}
T_{122}=1,\;\mbox{and}\;T_{ijk}=0\,\;\mbox{for other}\;i,j,k\in\{1,2\}.
\end{eqnarray*}
The eigenvalue equation of tensor $T$ becomes
\begin{eqnarray*}
\left\{\begin{array}{ccc}x_2^2&=&\lambda x_1^2,\\0&=&\lambda x_2^2.\end{array}\right.
\end{eqnarray*}
Obviously, $\rho(T)=0$. Hence, a nonzero nonnegative tensor may have zero spectral radius. So, Theorem \ref{add-thm-1} is not vacuous in general.
\end{Example}

At the end of this section, we summarize the Perron-Frobenius Theorem for nonnegative tensors as follows.

\begin{Theorem}\label{add-thm-2}
Let $T$ be an $m$-th order $n$ dimensional nonnegative tensor, then
\begin{itemize}
\item (Yang and Yang \cite{yy}) $\rho(T)$ is an eigenvalue of $T$ with a nonnegative eigenvector;
\item (Theorem \ref{add-thm-1}) if furthermore $T$ is strictly nonnegative, then $\rho(T)>0$;
\item (Friedland, Gaubert and Han \cite{fgh}) if furthermore $T$ is weakly irreducible, then $\rho(T)$ has a positive eigenvector;
\item (Chang, Pearson and Zhang \cite{cpz}) if furthermore $T$ is irreducible and if $\lambda$ is an eigenvalue with a nonnegative eigenvector, then $\lambda=\rho(T)$;
\item (Yang and Yang \cite{yy}) if $T$ is irreducible, and $T$ has $k$ distinct eigenvalues of modulus $\rho(T)$, then the eigenvalues are $\rho(T)\exp(i2\pi j/k)$ with $i^2=-1$ and $j=0,\ldots,k-1$;
\item (Chang, Pearson and Zhang \cite{cpz2}) if furthermore $T$ is primitive, then $k=1$; and
\item (Pearson \cite{p}) if $T$ is further essentially positive, $\rho(T)$ is real geometrically simple.
\end{itemize}
\end{Theorem}

\section{Relationships of the seven classes of nonnegative tensors}
\setcounter{Assumption}{0}
\setcounter{Theorem}{0} \setcounter{Proposition}{0}
\setcounter{Corollary}{0} \setcounter{Lemma}{0}
\setcounter{Definition}{0} \setcounter{Remark}{0}
\setcounter{Algorithm}{0}  \setcounter{Example}{0}
\hspace{4mm} In this section, we make a clear diagram of the
relationships among the concepts of nonnegative tensors mentioned above. Let $T$ be a
nonnegative tensor of order $m$ and dimension $n$ throughout this
section.  Let $E$ be the identity tensor of order $m$ and dimension
$n$ with its diagonal elements as $1$ and off-diagonal elements as
$0$ (when $m=2$, $E$ is the usual identity matrix).

The relationships among irreducibility, primitivity, weak positivity
and essential positivity were characterized in \cite[Section 3]{zqx}, which can
be summarized as follows:
\begin{itemize}
\item If $T$ is essentially positive, then $T$ is both weakly positive and
primitive. But the converse is not true. Moreover, there exists a tensor which is weakly positive and
primitive simultaneously, but not essentially positive. 
\item There is no inclusion relation between
the class of weakly positive tensors and the class of primitive
tensors.
\item If $T$ is weakly positive or primitive, then $T$ is
irreducible. But the converse is not true.
\end{itemize}

Actually, there are close relationships between (weakly)
irreducibility and (weakly) primitivity, and between weak positivity
and essential positivity.

\begin{Theorem}\label{thm-ad01}
A nonnegative tensor $T$ is weakly irreducible / irreducible /
weakly positive if and only if $T+E$ is weakly primitive / primitive
/ essentially positive.   A nonnegative tensor $T$ is essentially positive if and only if it is weakly positive, and all of its diagonal
elements are positive.
\end{Theorem}

\noindent {\bf Proof.}
\begin{itemize}
\item (Weakly irreducible / Weakly primitive) We have that nonnegative representation
matrix $G(T)$ is irreducible if and only if matrix $G(T+E)$ is primitive
\cite[Theorem 2.1.3 and Corollary 2.4.8]{bp}. So, by Definition \ref{def-new}, tensor $T$ is weakly
irreducible if and only if tensor $T+E$ is weakly primitive.
\item (Weakly positive / Essentially positive) The claims follow directly from Definition \ref{def-basic} and the nonnegativity of
$T$.
\item (Irreducible / Primitive) First, it follows from \cite[Theorem 6.6]{yy} immediately that $T$ is irreducible if and only if $F_{T+E}^{n-1}(x)>0$ for any nonzero $x\in\Re^n_+$.

    Second, if $T+E$ is primitive, then $T+E$ is irreducible. We can prove that
$F_{T+E}^{n-1}(x)>0$ for any nonzero $x\in\Re^n_+$, which implies that $T$ is irreducible by the above result. Actually, let $K:=T+\frac 1 2 E$, then, $2 K$ is irreducible and nonnegative. We have $2(T+E)=2(K+\frac 1 2 E)=2K+E$ and $F_{2 (T+E)}^{n-1}(x)=F_{2K+E}^{n-1}(x)>0$ for any nonzero $x\in\Re^n_+$ by the above result. It is straightforward to check that $F_{2 (T+E)}^{n-1}(x)=2^{\frac{n-1}{m-1}}F_{T+E}^{n-1}(x)$ for any nonzero $x\in\Re^n_+$. So, the result follows. \ep
\end{itemize}

We now discuss some other relations among
these six classes of nonnegative tensors.
By Definitions \ref{def-basic} and \ref{def-new}, if $T$ is
irreducible, then it is weakly irreducible. Nevertheless, the
converse is not true in general, which can be seen from the
following example.

\begin{Example}\label{exam4}
Let $T$ be a third order three dimensional tensor which is defined by
$T_{123}=T_{221}=T_{223}=T_{312}=T_{332}=1$ and $T_{ijk}=0$ for
other $i,j,k\in\{1,2,3\}$. Then,
$G(T)=\left(\begin{array}{ccc}0&1&1\\1&2&1\\1&2&1\end{array}\right)$
is irreducible but $T_{2ij}=0$ for all $i,j\in\{1,3\}$, which says
that $T$ is reducible.
\end{Example}

By Definition \ref{def-new}, if $T$ is weakly
primitive, then $T$ is weakly irreducible. The following example
demonstrates that the converse is not true.
\begin{Example}\label{exam4-a}
Let $T$ be a third order three dimensional tensor which is defined by $T_{122}=T_{233}=T_{311}=1$ and
$T_{ijk}=0$ for other $i,j,k\in\{1,2,3\}$. Then, $T$ is not  weakly primitive,
since its representation matrix
$G(T)=\left(\begin{array}{ccc}0&1&0\\0&0&1\\1&0&0\end{array}\right)$
is not primitive. But it is an irreducible matrix, and hence, $T$ is weakly
irreducible by Definition \ref{def-new}.
\end{Example}

We now discuss the relations between primitivity and weak
primitivity. The following result is a complementary to that in
Chang, Pearson and Zhang \cite{cpz2}.
\begin{Lemma}\label{lem-0}
For a nonnegative tensor $T$ of order $m$ and dimension $n$, if
$M(T)$ is primitive, then $T$ is primitive, and if $T$ is primitive, then
$G(T)$ is primitive.
\end{Lemma}

\noindent {\bf Proof.} If $M(T)$ is primitive,  then let
$K:=\left[M(T)\right]^{k}$ with $k:=n^2-2n+2$, we have $K_{ij}>0$
for any $i,j\in\{1,\ldots,n\}$ \cite[Theorem 2.4.14]{bp}. Now, for any nonzero
$x\in\Re^n_+$, suppose $x_j>0$. Then, for any $i\in\{1,\ldots,n\}$,
there exist $i_2,\ldots,i_{k}$ such that $M(T)_{ii_2},
M(T)_{i_2i_3},\ldots,M(T)_{i_{k}j}>0$. So, $T_{ii_2\ldots
i_2},\ldots,T_{i_{k}j\ldots j}>0$. Thus, we have
$\left[F_T^{k}(x)\right]_i>0$ for any $i\in\{1,\ldots,n\}$. As the
above inequalities hold for any nonzero $x\in\Re^n_+$, we obtain
that $T$ is primitive.

If $T$ is primitive, then for some integer $k>0$, $F_T^{k}(x)>0$ for any nonzero $x\in\Re^n_+$. For any $i\in\{1,\ldots,n\}$, let $e_j$ denote
the $j$-th column of the $n\times n$ identity matrix for any
$j\in\{1,\ldots,n\}$. We thus have $\left[F_T^{k}(e_j)\right]_i>0$.
So, we must have indices $\{i_2^2,\ldots,i_m^2\}$,
$\{i_1^3,\ldots,i_m^3\}$, $\ldots$,
$\{i_1^{k-1},\ldots,i_m^{k-1}\}$, $i_k$ such that $T_{ii_2^2\ldots
i_m^2},T_{i_1^{3}\ldots i_m^{3}},\ldots,T_{i_1^{k-1}\ldots
i_m^{k-1}},T_{i_k j\ldots j}>0$ and $i^{l+1}\in\{i^l_2,\ldots,
i^l_m\}$ for $l\in\{1,\ldots,k-1\}$. Thus, if we let
$L:=\left[G(T)\right]^{k}$, we should have $L_{ij}>0$ for all
$i,j\in\{1,\ldots,n\}$. Hence, $G(T)$ is primitive. \ep

By Lemma \ref{lem-0} and Definition \ref{def-new}, if $T$ is
primitive, then $T$ is weakly primitive.

We list some further relationships among the concepts mentioned
above as follows:
\begin{itemize}
\item By Example \ref{exam4}, we see that $G(T)$ is primitive but $T$ is
reducible since $T_{2ij}=0$ for all $i,j\in\{1,3\}$. Hence, there
exists a nonnegative tensor which is weakly primitive but not
irreducible.
\item Let $T$ be a third order two dimensional tensor which is defined by $T_{122}=T_{211}=1$ and $T_{ijk}=0$ for other
$i,j,k\in\{1,2\}$. Then, $T$ is weakly positive, and hence,
irreducible, but not weakly primitive.
\item Let $T$ be a third order two dimensional tensor which is defined by
$T_{122}=T_{211}=T_{212}=T_{121}=1$ and $T_{ijk}=0$ for other
$i,j,k\in\{1,2\}$. Then, $T$ is weakly positive but not primitive,
since $Te_1^2=e_2$ and $Te_2^2=e_1$. However, it is weakly primitive,
since $G(T)=\left(\begin{array}{cc}1&2\\2&1\end{array}\right)$ is
obviously primitive.
\end{itemize}

Using Corollary \ref{add-cor-1} and Example \ref{add-exm-1}, we
can summarize the relationships obtained so far in Figure 1.

\begin{figure}[htbp]
\centering
\includegraphics[width=5.8in]{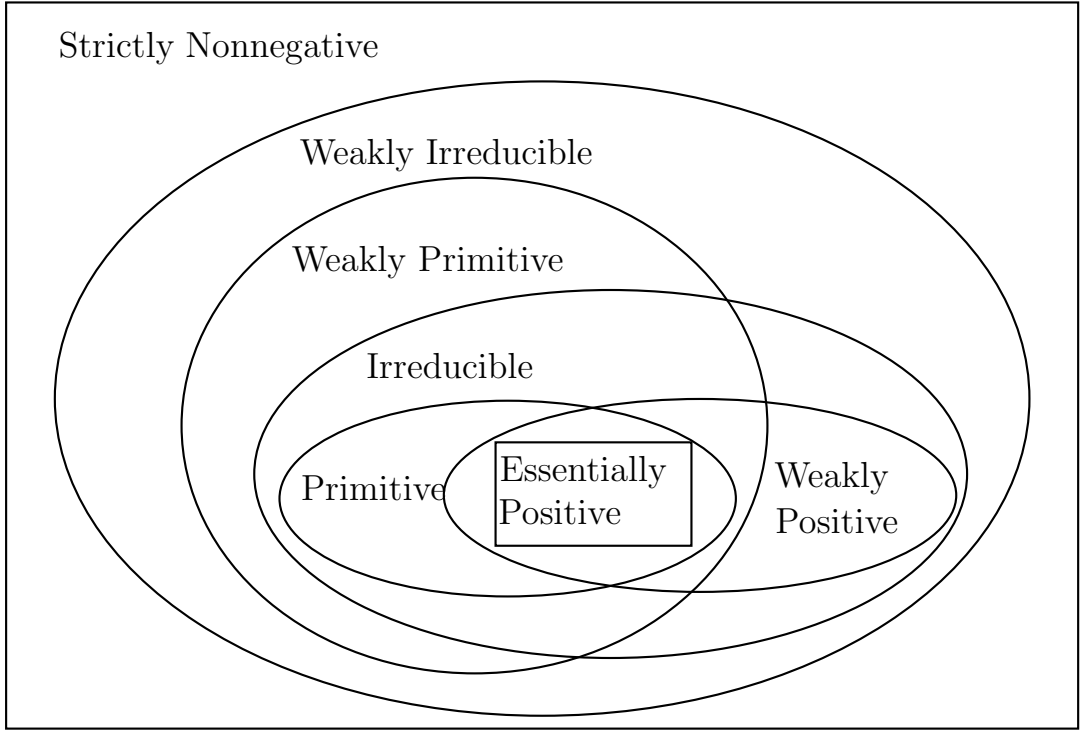}
\caption{Relationships of the seven classes of nonnegative tensors}
\end{figure}

\section{Global R-linear convergence of a power method for weakly irreducible nonnegative tensors}
\setcounter{Assumption}{0}
\setcounter{Theorem}{0} \setcounter{Proposition}{0}
\setcounter{Corollary}{0} \setcounter{Lemma}{0}
\setcounter{Definition}{0} \setcounter{Remark}{0}
\setcounter{Algorithm}{0}  \setcounter{Example}{0}
We present here a modification of the power method proposed in
\cite{fgh}.

\hspace{4mm}

\begin{Algorithm}\label{algo} (A Higher Order Power Method (HOPM))
\begin{description}
\item [Step 0] Initialization: choose $x^{(0)}\in \Re^n_{++}$.  Let $k:=0$.

\item [Step 1] Compute
\begin{eqnarray*}
\begin{array}{c}
\bar x^{(k+1)}:=T(x^{(k)})^{m-1},\quad x^{(k+1)}:=\frac{\left(\bar
x^{(k+1)}\right)^{[\frac{1}{m-1}]}}{e^T\left[\left(\bar
x^{(k+1)}\right)^{[\frac{1}{m-1}]}\right]},\\
\alpha\left(x^{(k+1)}\right):=\max_{1\leq i\leq
n}\frac{\left(T(x^{(k)})^{m-1}\right)_i}{\left(x^{(k)}\right)_i^{m-1}}\quad
\mbox{\rm and}\quad \beta\left(x^{(k+1)}\right):=\min_{1\leq i\leq
n}\frac{\left(T(x^{(k)})^{m-1}\right)_i}{\left(x^{(k)}\right)_i^{m-1}}.
\end{array}
\end{eqnarray*}

\item [Step 2] If $\alpha\left(x^{(k+1)}\right)=\beta\left(x^{(k+1)}\right)$, stop. Otherwise, let $k:=k+1$, go to Step 1.
\end{description}
\end{Algorithm}

Algorithm \ref{algo} is well-defined if the underlying tensor $T$ is
a strictly nonnegative tensor, as in this case, $Tx^{m-1}>0$ for any
$x>0$.  Hence, Algorithm \ref{algo} is also well-defined for weakly
irreducible nonnegative tensors. The following theorem establishes
convergence of Algorithm \ref{algo} if the underlying tensor $T$ is
weakly primitive, where we need to use the concept of Hilbert's
projective metric \cite{n}. We first recall such a concept. For any
$x,y\in\Re^n_{+}\setminus\{0\}$, if there are $\alpha,\beta>0$ such
that $\alpha x\leq y\leq\beta x$, then $x$ and $y$ are called {\em
comparable}. If $x$ and $y$ are comparable, and define
\begin{eqnarray*}
m(y/x):=\sup\{\alpha>0\;|\;\alpha x\leq y\}\quad\mbox{and}\quad
M(y/x):=\inf\{\beta>0\;|\;y\leq \beta x\},
\end{eqnarray*}
then, the Hilbert's projective metric $d$ can be defined by
\begin{eqnarray*}
d(x,y):=\left\{\begin{array}{ll} \mbox{log}(\frac{M(y/x)}{m(y/x)}),
&\mbox{if}\; x \;\mbox{and}\; y \;\mbox{are}\;
\mbox{comparable},\\+\infty, &\mbox{otherwise}\end{array}\right.
\end{eqnarray*}
for $x,y\in\Re^n_{+}\setminus\{0\}$. Note that if
$x,y\in\Delta_n:=\{z\in\Re^n_{++}\;|\;e^Tz=1\}$, then $d(x,y)=0$ if and
only if $x=y$. Actually, it is easy to check that $d$ is a metric on
$\Delta_n$.

\begin{Theorem}\label{thm-ad0}
Suppose that $T$ is a weakly irreducible nonnegative tensor of order
$m$ and dimension $n$.  Then, the following results hold.
\begin{itemize}
\item[(i)] $T$ has a positive eigenpair $(\lambda,x)$, and $x$ is unique up
to a multiplicative constant.
\item[(ii)] Let $(\lambda_*,x^*)$ be the unique positive
eigenpair of $T$ with $\sum_{i=1}^n(x^*)_i=1$. Then,
$$
\min_{x\in\Re^n_{++}}\max_{1\leq i\leq
n}\frac{\left(Tx^{m-1}\right)_i}{x_i^{m-1}}=\lambda_*=\max_{x\in\Re^n_{++}}\min_{1\leq
i\leq n}\frac{\left(Tx^{m-1}\right)_i}{x_i^{m-1}}.
$$
\item[(iii)] If $(\nu, v)$ is another eigenpair of $T$, then $|\nu|\leq
\lambda_*$.
\item[(iv)] Suppose that $T$ is weakly primitive and the sequence $\{x^{(k)}\}$ is generated by
Algorithm \ref{algo}. Then, $\{x^{(k)}\}$ converges to the unique vector
$x^*\in \Re^n_{++}$ satisfying $T(x^*)^{m-1}=\lambda_*
(x^*)^{[m-1]}$ and $\sum_{i=1}^nx_i^*=1$, and there exist constant $\theta\in(0,1)$ and positive integer $M$ such that
\begin{eqnarray}\label{grline}
d(x^{(k)},x^*)\leq\theta^{\frac{k}{M}} \frac{d(x^{(0)},x^*)}{\theta}
\end{eqnarray}
holds for all $k\geq 1$.
\end{itemize}
\end{Theorem}

\noindent {\bf Proof.} Except the result in \reff{grline}, all other
results in this theorem can be easily obtained from \cite[Theorem
4.1, Corollaries 4.2, 4.3 and 5.1]{fgh}. So, we only give the proof of
\reff{grline} here. We have the following observations first:
\begin{itemize}
\item $\Re^n_{+}$ is a {\it normal cone} in Banach space $\Re^n$, since $y\geq x\geq 0$ implies $\|y\|\geq\|x\|$;
\item $\Re^n_+$ has nonempty interior $\Re^n_{++}$ which is an open cone,
and $F_T:\Re^n_{++}\ra\Re^n_{++}$ is continuous and {\it
order-preserving} by Corollary \ref{add-cor-1} and the nonnegativity of
tensor $T$;
\item $F_T$ is homogeneous of degree $1$ in $\Re^n_{++}$;
\item the set $\Delta_n$ is connected and $T$ has an eigenvector $x^*$ in $\Delta_n$ by Theorem \ref{thm-ad0} (i);
\item by \reff{map-1}, $F_T$ is continuously differentiable in an open neighborhood of $x^*$, since $x^*>0$;
\item by Definition \ref{def-new}, $G(T)$ is primitive, hence there exists an integer $N$ such that $\left[G(T)\right]^N>0$. So, $\left[G(T)\right]^Nx$ is comparable with $x^*$ for any nonzero $x\in\Re^n_+$;
\item $G(T):\Re^n\ra\Re^n$ is a compact linear map, hence its essential spectrum radius is zero \cite[Page 38]{n}, while its spectral radius is positive since it is a primitive matrix \cite{bp}.
\end{itemize}
Hence, by \cite[Corollary 2.5 and Theorem 2.7]{n}, we have that
there exist a constant $\theta\in(0,1)$ and a positive integer $M$
such that
\begin{eqnarray}\label{grline1}
d(x^{(Mj)},x^*)\leq\theta^j d(x^{(0)},x^*),
\end{eqnarray}
where $d$ denotes the Hilbert's projective metric on
$\Re^n_+\setminus\{0\}$.

By \cite[Proposition 1.5]{n}, we also have that
\begin{eqnarray}\label{contr}
d(F_T(x),F_T(y))\leq d(x,y)
\end{eqnarray}
for any $x,y\in\Re^n_+$. Since $\lambda_*>0$, by the property of
Hilbert's projective metric $\mbox{d}$ \cite[Page 13]{n}  we have
that
\begin{eqnarray*}
d(x^{(k+1)},x^*)&=&d\left(\frac{F_T(x^{(k)})}{e^TF_T(x^{(k)})},x^*\right)
=d\left(\frac{F_T(x^{(k)})}{e^TF_T(x^{(k)})},\frac{1}{(\lambda_*)^{\frac{1}{m-1}}}F_T(x^*)\right)\\
&=&d(F_T(x^{(k)}),F_T(x^*))\leq d(x^{(k)},x^*)
\end{eqnarray*}
holds for any $k$. So, for any $k\geq M$, we could find the largest
$j$ such that $k\geq Mj$ and $M(j+1)\geq k$. Hence,
\begin{eqnarray*}
d(x^{(k)},x^*)\leq d(x^{(Mj)},x^*)\leq \theta^{j}d(x^{(0)},x^*)\leq \theta^{\frac{k}{M}-1}d(x^{(0)},x^*)
\end{eqnarray*}
which implies \reff{grline} for all $k\geq M$. When $1\leq k<M$, we
have $\theta^{\frac{k}{M}}>\theta$, since $\theta\in(0,1)$.
Therefore, \reff{grline} is true for all $k\geq 1$. \ep

We denote by $x^{[p]}$ a vector with its $i$-th element being
$x_i^p$.

By Theorems \ref{thm-ad01} and \ref{thm-ad0}, the following result
holds obviously.
\begin{Theorem}\label{thm}
Suppose that $T$ is a weakly irreducible nonnegative tensor of order
$m$ and dimension $n$, and the sequence $\{x^{(k)}\}$ is generated by
Algorithm \ref{algo} with $T$ being replaced by $T+E$. Then,
$\{x^{(k)}\}$ converges to the unique vector $x^*\in \Re^n_{++}$
satisfying $T(x^*)^{m-1}=\lambda_* (x^*)^{[m-1]}$ and
$\sum_{i=1}^nx_i^*=1$, and there exist a constant $\theta\in(0,1)$
and a positive integer $M$ such that \reff{grline} holds for all
$k\geq 1$.
\end{Theorem}

\begin{Remark}
(i) Compared with \cite[Corollaries 5.1 and 5.2]{fgh}, a main
advantage of our results is that \reff{grline} in Theorem
\ref{thm-ad0}(iv) gives the {\bf global} $R$-linear convergence of
Algorithm \ref{algo}; while the geometric convergence given in
\cite[Corollary 5.2]{fgh} is essentially a result of {\bf local}
$R$-linear convergence. (ii) Compared with the results in
\cite{zqx}, a main advantage of our results is that the problem we
considered in Theorem \ref{thm} is broader than that in \cite{zqx},
i.e., our results are obtained for the weakly irreducible
nonnegative tensors; while the results in \cite{zqx} hold for the
irreducible nonnegative tensors. There are also other differences
between Theorem \ref{thm} and those in \cite{zqx}, such as, Theorem
\ref{thm} gives the global $R$-linear convergence of the {\bf
iterated sequence $\{x^{(k)}\}$} of Algorithm \ref{algo} for weakly irreducible nonnegative tenors; while the
Q-linear convergence of the corresponding {\bf eigenvalue sequence}
was proved in \cite{zqx} for weakly positive nonnegative tensors.
\end{Remark}

Theorem \ref{thm} is an important basis for us to develop a method
for finding the spectral radius of a {\it general nonnegative
tensor}.

\section{Partition a general nonnegative tensor to weakly irreducible nonnegative tensors}
\setcounter{Assumption}{0}
\setcounter{Theorem}{0} \setcounter{Proposition}{0}
\setcounter{Corollary}{0} \setcounter{Lemma}{0}
\setcounter{Definition}{0} \setcounter{Remark}{0}
\setcounter{Algorithm}{0}  \setcounter{Example}{0}
\hspace{4mm} If a nonnegative tensor $T$ of order $m$ and dimension
$n$ is weakly irreducible, then from Theorem \ref{thm}, we can
find the spectral radius and the corresponding positive eigenvector
of $T$ by using Algorithm \ref{algo}. A natural question is that,
{\it if $T$ is not weakly irreducible, what can we do for it}?

In this section, we show that if a nonnegative tensor $T$ is not
weakly irreducible, then there exists a partition of the index set
$\{1,\ldots,n\}$ such that every tensor induced by the set in the
partition is weakly irreducible; and the largest eigenvalue of $T$
can be obtained from these induced tensors. Thus, we can find the
spectral radius of a general nonnegative tensor by using Algorithm
\ref{algo} for these induced weakly irreducible tensors. At the end
of this section, we show that, if weakly irreducibility is replaced
by irreducibility, a similar method does not work.

The following result is the theoretical basis of our method.

\begin{Theorem}\cite[Theorem 2.3]{yy}\label{thm-yy}
For any nonnegative tensor $T$ of order $m$ and dimension $n$,
$\rho(T)$ is an eigenvalue with a nonnegative eigenvector
$x\in\Re^n_+$ corresponding to it.
\end{Theorem}

To develop an algorithm for general nonnegative tensors, we prove
the following theorem which is an extension of the corresponding
result for nonnegative matrices \cite{bp}. For the convenience of the sequel analysis, a one dimensional tensor is always considered as irreducible, hence weakly irreducible. Similarly, one dimensional positive tensors are considered as primitive. Note that Algorithm \ref{algo} works for one dimensional primitive tensor as well. Nonetheless, weakly irreducible nonnegative tensors with dimension one may have zero spectral radius, but they are always positive when the dimension $n\geq 2$ by Theorem \ref{add-thm-1}. Note that, $n$ is assumed to be no smaller than two throughout this paper, while the case of one dimensional tensors is needed in the presentation of partition results in this section.

\begin{Theorem}\label{thm-5}
Suppose that $T$ is a nonnegative tensor of order $m$ and dimension
$n$. If $T$ is weakly reducible, then there is a partition
$\{I_1,\ldots,I_k\}$ of $\{1,\ldots,n\}$ such that every tensor in
$\{T_{I_j}\;|\;j\in\{1,\ldots,k\}\}$ is weakly irreducible.
\end{Theorem}

\noindent {\bf Proof.} Since $T$ is weakly reducible, by Definition
\ref{def-new} we can obtain that the matrix $G(T)$ is reducible.
Thus, we could find a partition $\{J_1,\ldots,J_l\}$ of
$\{1,\ldots,n\}$ such that
\begin{itemize}
\item [($\star$)] every matrix (a second order tensor) in
$\{\left[G(T)\right]_{J_i}\;|\;i\in\{1,\ldots,l\}\}$ is irreducible
and $\left[G(T)\right]_{st}=0$ for any $s\in J_p$ and $t\in J_q$
such that $p>q$.
\end{itemize}
Actually, by the definition of reducibility of a matrix, we can
find a partition $\{J_1,J_2\}$ of $\{1,\ldots,n\}$ such that
$\left[G(T)\right]_{st}=0$ for any $s\in J_2$ and $t\in J_1$. If
both $\left[G(T)\right]_{J_1}$ and $\left[G(T)\right]_{J_2}$ are
irreducible, then we are done. Otherwise, we can repeat the above
analysis to any reducible block(s) obtained above. In this way,
since $\{1,\ldots,n\}$ is a finite set, we can arrive at the
desired result ($\star$).

Now, if every tensor in $\{T_{J_i}\;|\;i\in\{1,\ldots,l\}\}$ is
weakly irreducible, then we are done. Otherwise, we repeat the above
procedure to generate a partition of $T$ to these induced tensors
which are not weakly irreducible. Since $\{1,\ldots,n\}$ is finite,
this process will stop in finite steps. Hence, the theorem
follows. \ep

By Theorems \ref{thm-5} and \ref{thm-ad01}, we have the following
corollary.
\begin{Corollary}\label{cor-2}
Suppose that $T$ is a nonnegative tensor of order $m$ and dimension
$n$. If $T$ is weakly irreducible, then $T+E$ is weakly primitive;
otherwise, there is a partition $\{I_1,\ldots,I_k\}$ of
$\{1,\ldots,n\}$ such that every tensor in
$\{(T+E)_{I_j}\;|\;j\in\{1,\ldots,k\}\}$ is weakly primitive.
\end{Corollary}

Given a nonempty subset $I$ of $\{1,\ldots,n\}$ and an $n$ vector $x$, we
will denote by $x_I$ an $n$ vector with its $i$-th element being
$x_i$ if $i\in I$ and zero otherwise; and $x(I)$ a $|I|$ vector
after deleting $x_j$ for $j\notin I$ from $x$.

\begin{Theorem}\label{thm-6}
Suppose that $T$ is a weakly reducible nonnegative tensor of order
$m$ and dimension $n$, and $\{I_1,\ldots,I_k\}$ is the partition of
$\{1,\ldots,n\}$ determined by Theorem \ref{thm-5}. Then,
$\rho(T)=\rho(T_{I_p})$ for some $p\in\{1,\ldots,k\}$.
\end{Theorem}

\noindent {\bf Proof.} By the proof of Theorem \ref{thm-5}, for the
nonnegative matrix $G(T)$, we could find a partition
$\{J_1,\ldots,J_l\}$ of $\{1,\ldots,n\}$ such that
\begin{itemize}
\item every matrix in
$\{\left[G(T)\right]_{J_i}\;|\;i\in\{1,\ldots,l\}\}$ is irreducible
and $\left[G(T)\right]_{st}=0$ for any $s\in J_p$ and $t\in J_q$
such that $p>q$.
\end{itemize}

First, we have that $\rho(T_{J_i})\leq\rho(T)$ for all $i\in\{1,\ldots,l\}$ by Lemma \ref{add-lem-2}.

Then, denote by $(\rho(T),x)$ a nonnegative
eigenpair of $T$ which is guaranteed by Theorem \ref{thm-yy}. Since
$\left[G(T)\right]_{ij}=0$ for all $i\in J_l$ and
$j\in\cup_{s=1}^{l-1}J_s$. We must have
\begin{eqnarray}\label{rela-2}
T_{ii_2\ldots i_m}=0\;\;\forall i\in J_l,\;\;\forall
\{i_2,\ldots,i_m\}\not\subseteq J_l.
\end{eqnarray}
Hence, for all $i\in J_l$, we have
\begin{eqnarray*}
\rho(T)x_i^{m-1}&=&(Tx^{m-1})_i\\
&=&\sum_{i_2,\ldots,i_m=1}^nT_{ii_2\ldots i_m}x_{i_2}\cdots x_{i_m}\\
&=&\sum_{\{i_2,\ldots,i_m\}\subseteq J_l}^nT_{ii_2\ldots i_m}x_{i_2}\cdots x_{i_m}\\
&=&\left\{T_{J_l}\left(x(J_l)\right)^{m-1}\right\}_i,
\end{eqnarray*}
where the third equality follows from \reff{rela-2}. If $x(J_l)\neq
0$, then $(\rho(T),x(J_l))$ is a nonnegative eigenpair of tensor
$T_{J_l}$; and if $x(J_l)=0$, then we have
\begin{eqnarray*}
T_{\cup_{j=1}^{l-1}J_j}\left(x({\cup_{j=1}^{l-1}J_j})\right)^{m-1}=\rho(T)\left[x({\cup_{j=1}^{l-1}J_j})\right]^{[m-1]}.
\end{eqnarray*}
In the later case, repeat the above analysis with $T$ being replaced
by $T_{\cup_{j=1}^{l-1}J_j}$. Since $x\neq 0$ and $l$ is finite, we
must find some $t\in\{1,\ldots,l\}$ such that $x(J_t)\neq 0$ and
$(\rho(T),x(J_t))$ is a nonnegative eigenpair of tensor $T_{J_t}$.

Now, if $T_{J_t}$ is weakly irreducible, we are done since $J_t=I_p$
for some $p\in\{1,\ldots,k\}$ by the proof of Theorem \ref{thm-5}.
Otherwise, repeat the above analysis with $T$ and $x$ being replaced
by $T_{J_t}$ and $x(J_t)$, respectively. Such a process is finite,
since $n$ is finite. Thus, we always obtain a weakly irreducible
nonnegative tensor $T_{I_p}$ with $I_p\subseteq \{1,\ldots,n\}$ for
some $p$ such that $(\rho(T),x(I_p))$ is a nonnegative eigenpair of
tensor $T_{I_p}$. Furthermore, $(\rho(T),x_{S})$ with $S:=\cup_{i=1}^pI_p$ is a nonnegative
eigenpair of tensor $T$.

The proof is complete. \ep

Note that if $T$ is furthermore symmetric, then we could get a diagonal block representation of $T$ with diagonal blocks $T_{I_i}$ (after some permutation, if necessary).
Now, by Corollary \ref{cor-2} and Theorems \ref{thm}, \ref{thm-5}
and \ref{thm-6}, we could get the following theorem.

\begin{Theorem}\label{thm-7}
Suppose that $T$ is a nonnegative tensor of order $m$ and dimension
$n$.
\begin{itemize}
\item [(a)] If $T$ is weakly irreducible, then $T+E$ is weakly primitive by Theorem \ref{thm-ad01};
and hence, Algorithm \ref{algo} with $T$ being replaced by $T+E$
converges to the unique positive eigenpair $(\rho(T+E),x)$ of $T+E$.
Moreover, $(\rho(T+E)-1,x)$ is the unique positive eigenpair of $T$.
\item [(b)] If $T$ is not weakly irreducible, then, we can get a set of
weakly irreducible tensors $\{T_{I_j}\;|\;j=1,\ldots,k\}$ with $k>1$
by Theorem \ref{thm-5}. For each $j\in\{1,\ldots,k\}$, we use item
(a) to find the unique positive eigenpair $(\rho(T_{I_j}),x^{j})$ of
$T_{I_j}$ which is guaranteed by Corollary \ref{cor-2} when $|I_j|\geq 2$ or eigenpair $(T_{I_j},1)$ when $|I_j|=1$. Then,
$\rho(T)=\max_{j=1,\ldots,k}\rho(T_{I_j})$ by Theorem
\ref{thm-6}. If $T$ is furthermore symmetric, then, $x$ with
$x(I_{j^*})=x^{{j^*}}$ is a nonnegative eigenvector of $T$ where
$j^*\in\mbox{argmax}_{j=1,\ldots,k}\rho(T_{I_j})$
\end{itemize}
\end{Theorem}

\begin{Remark}\label{rmk-2}
By Theorem \ref{thm-7} and Algorithm \ref{algo}, if nonnegative
tensor $T$ is weakly irreducible, then, the spectral radius of $T$
can be found directly by Algorithm \ref{algo} with $T$ being
replaced by $T+E$. If $T$ is not weakly irreducible, then, we have to
find the partition of $\{1,\ldots,n\}$ determined by Theorem
\ref{thm-5}. Fortunately, we can find such a partition through the
corresponding partition of the nonnegative representation matrix of
$T$ and its induced tensors according to Theorem \ref{thm-6}. The
specific method of finding such a partition is given in the next
section.
\end{Remark}

Most of the known papers, which established the Perron-Frobenius
theorem and showed the convergence of the power method for
nonnegative tensors, concentrated on the class of irreducible
nonnegative tensors. A natural question is that, {\em for any given
reducible nonnegative tensor $T$, whether a partition of $T$ similar
to the result given in Theorem \ref{thm-5} can be derived or not}.
If so, {\em whether the spectral radius of $T$ can be obtained by
using the power method for the induced irreducible nonnegative
tensors or not}. At the end of this section, we answer these two questions.
The answer to the first question is positive, which is given as
follows.
\begin{Theorem}\label{thm-4}
Suppose that $T$ is a nonnegative tensor of order $m$ and dimension
$n$. If $T$ is reducible, then there is a partition
$\{I_1,\ldots,I_k\}$ of $\{1,2,\ldots,n\}$ such that any one of the
tensors $\{T_{I_j}\;|\;j\in\{1,\ldots,k\}\}$ is irreducible and
\begin{eqnarray*}
T_{st_2\ldots t_m}=0,\;\;\forall s\in I_p,\;\forall
\{t_2,\ldots,t_m\}\subset I_q,\;\forall p>q.
\end{eqnarray*}
\end{Theorem}

\noindent {\bf Proof.} Since $T$ is reducible, by the definition of
reducibility, there exists a nonempty proper subset $I_2$ of
$\{1,\ldots,n\}$ such that
\begin{eqnarray*}
T_{ii_2\ldots i_m}=0,\;\forall i\in I_2,\;\forall i_2,\ldots,i_m\in
I_1:=\{1,\ldots,n\}\setminus I_2.
\end{eqnarray*}
If both $T_{I_1}$ and $T_{I_2}$ are irreducible, then we are done.
Without loss of generality, we assume that $T_{I_1}$ is irreducible
and $T_{I_2}$ is reducible. Then, by the reducibility of $T_{I_2}$,
we can get a partition $\{J_2,J_3\}$ of $I_2$ such that
\begin{eqnarray*}
T_{ii_2\ldots i_m}=0,\;\forall i\in J_3,\;\forall i_2,\ldots,i_m\in
J_2:=I_2\setminus J_3.
\end{eqnarray*}
If both $T_{J_2}$ and $T_{J_3}$ are irreducible, then we are done,
since $\{I_1,J_2,J_3\}$ is the desired partition of $\{1,\ldots,n\}$.
Otherwise, repeating the above procedure, we can get the desired
results, since $n$ is finite.  \ep

However, the answer to the second question is negative, which can be
seen by the following example.
\begin{Example}\label{exm-2}
Let $T$ be a third order two dimensional tensor which is defined by
\begin{eqnarray*}
T_{111}=1,\;T_{112}=T_{121}=T_{211}=4,\;T_{122}=T_{212}=T_{221}=0,\;\mbox{and}\;T_{222}=1.
\end{eqnarray*}
Since $T_{122}=0$, tensor $T$ is reducible. And $T_{111}=T_{222}=1$
is the largest eigenvalue of both induced tensors by Theorem \ref{thm-4}. While the
nonnegative eigenpairs of $T$ are
\begin{eqnarray*}
(1,(0,1)^T)\quad \mbox{\rm and}\quad (7.3496, (0.5575,0.4425)^T).
\end{eqnarray*}
\end{Example}

This example prevents us to use $\rho(T_{I_i})$s' to get $\rho(T)$
under the framework of irreducibility. In addition, it is easy to
see that checking weak irreducibility of tensor $T$ is much easier
than checking irreducibility of tensor $T$, since the former is based
on a nonnegative matrix which has both sophisticated theory and
algorithms \cite{bp}.

\section{A specific algorithm and numerical experiments}
\setcounter{Assumption}{0}
\setcounter{Theorem}{0} \setcounter{Proposition}{0}
\setcounter{Corollary}{0} \setcounter{Lemma}{0}
\setcounter{Definition}{0} \setcounter{Remark}{0}
\setcounter{Algorithm}{0}  \setcounter{Example}{0}
\hspace{4mm} In this section, based on Algorithm \ref{algo}, Theorem
\ref{thm}, and the theory established in Section 5, we develop a
specific algorithm for finding the spectral radius of a general nonnegative
tensor.

\subsection{A specific algorithm}

\hspace{4mm} In this subsection, we give an algorithm for finding irreducible
blocks of a nonnegative matrix $M$, which is based on the fact that a
nonnegative matrix $M$ is irreducible if and only if
$\left(M+E\right)^{n-1}>0$ \cite{bp}.

\begin{Algorithm}\label{algo2}(Irreducible blocks of nonnegative matrices)
\begin{description}
\item [Step 0] Given a nonnegative matrix $M$, let $k=1$ and $C^1:=M+E$.

\item [Step 1] Until $k=n-1$, repeat $C^k:=C^{k-1}(M+E)$ and $k:=k+1$.

\item [Step 2] Sort the numbers of nonzero elements of columns of $C^{n-1}$
in ascend order, then perform symmetric permutation to $C^{n-1}$
according to the sorting order into a matrix $K$; sort the numbers
of nonzero elements of rows of $K$ in descend order, then perform
symmetric permutation to $K$ according to the sorting order into a
matrix $L$. Record the two sorting orders.
\item [Step 3] Let $i=1$, $j=1$ and $s=1$, create an index set $I_j$ and an vector $ind$,
put $i$ into $I_j$ and index it to be the $s$-th element in $I_j$,
and set $ind(i)$ to be $1$.
\item [Step 4] If $L(d,I_j(s))>0$ and $L(I_j(s),d)>0$ for some $d\in\{1,\ldots,n\}$
with $ind(d)$ being not $1$, where $I_j(s)$ is the $s$-th element in
$I_j$, set $s:=s+1$ and put $d$ into $I_j$ and index it to be the
$s$-th element in $I_j$. Set $ind(d)$ to be $1$, and $i:=i+1$. If
there is no such $d$ or $i=n$, go to Step 6; if there is no such $d$
but $i<n$, go to Step 5.
\item [Step 5] Let $j:=j+1$ and $s:=1$. Create index set $I_j$,
and find a $d$ with $ind(d)$ being not $1$, put $d$ into $I_j$ and
index it to be the $s$-th element in $I_j$. Set $ind(d)$ to be $1$,
go back to Step 4.
\item [Step 6] Using sorting orders in Step 2 and partition $\{I_1,\ldots,I_k\}$
found by Steps 3-5, we could find the partition for the matrix $M$ easily.
\end{description}
\end{Algorithm}

Now, we propose a specific algorithm for finding the spectral radius of a general
nonnegative tensor.

\begin{Algorithm}\label{algo3}(Spectral radius of a nonnegative tensor)
\begin{description}
\item [Step 0] Let $v$ be an $n$ vector with its elements being zeros, and $i=1$.

\item [Step 1] Given a nonnegative tensor $T$ of order $m$ and dimension $n$,
compute its representation matrix $G(T)$ as Definition
\ref{def-new}.

\item [Step 2] Finding out the partition $\{I_1,\ldots,I_k\}$ of $\{1,\ldots,n\}$
using Algorithm \ref{algo2} with $M$ being replaced by $G(T)$.

\item [Step 3] If $k=1$, using Algorithm \ref{algo} to finding out the spectral
radius $\rho$ of $T$,
set $v(i)=\rho$, and $i=i+1$.
Otherwise, go to Step 4.

\item [Step 4] For $j=1,\ldots,k$, computing the induced tensor $T_{I_j}$
and its corresponding representation matrix $G_j$, set $T$ as
$T_{I_j}$, $G(T)$ as $G_j$ and $n$ as $|I_j|$, run subroutine
Steps 2-4.

\item [Step 5] Out put the spectral radius of $T$ as $\max_{i=1}^nv(i)$.
\end{description}
\end{Algorithm}

\subsection{Numerical experiments}

\hspace{4mm} In this subsection, we report some preliminary numerical
results for computing the spectral radius of a general nonnegative tensor using Algorithm
\ref{algo3} (with Algorithms \ref{algo} and \ref{algo2}). All
experiments are done on a PC with CPU of 3.4 GHz and RAM of 2.0 GB,
and all codes are written in MATLAB.

To demonstrate that Algorithm \ref{algo3} works for general
nonnegative tensors, we randomly generate third order nonnegative
tensors with dimension $n$ which is specialized in Table 1. We
generate the testing tensors by randomly generating their every element
uniformly in $[0,1]$ with a density {\bf Den} which is specialized
in Table 1. We use Algorithm \ref{algo3} to find the spectral radii of the generated tensors for every
case, i.e., with different dimensions $n$ and element density {\bf
Den}. The algorithm is terminated if $|\alpha(x^{(k)})-\beta(x^{(k)})|\leq
10^{-6}$. For every case, we simulate $50$ times to get the average
spectral radius {\bf $\rho:=\frac{\alpha(x^{(k)})+\beta(x^{(k)})}{2}$}, the
average number of iterations {\bf Ite} performed by Algorithm
\ref{algo}, the average weakly irreducible blocks of the generated
tensors {\bf Blks}, and the average residual of $Tx^2-\rho x^{[2]}$ with
the found spectral radius $\rho$ and its corresponding eigenvector
$x$ in $2$-norm {\bf Res}. We also use {\bf Per} to denote the
percentage of weakly irreducible tensors generated among the $50$
simulations, and {\bf TolCpu} to denote the total cputime spent for
the simulation in every case. All results are listed in Table 1. In
addition, we also test the following example and the related
numerical results are listed in Table 2, where {\bf Blk} denotes the
block number of the computation, {\bf Ite} denotes the iteration
number, and the other items are clear from Algorithm \ref{algo3}.
The initial points for both tests in Tables 1 and 2 are randomly
generated with their every element uniformly in $(0,1)$.
All the simulated tensors are strictly nonnegative by the above simulation
strategy, since every component of $R(T)$ is the summation of so many terms, it is never zero in our simulation.

\begin{Example}\label{exm-6}
Let third order three dimensional tensor $T$ be defined by $T_{111}=T_{222}=1$, $T_{122}=3$,
$T_{211}=5$, $T_{333}=4$ and $T_{ijk}=0$ for other
$i,j,k\in\{1,2,3\}$. Then, the eigenvalue problem \reff{eig} reduces to:
\begin{eqnarray*}
\left\{\begin{array}{rcl}x_1^2+3x_2^2&=&\lambda
x_1^2,\\5x_1^2+x_2^2&=&\lambda x_2^2,\\4x_3^2&=&\lambda
x_3^2.\end{array}\right.
\end{eqnarray*}
It is easy to see that tensor $T$ is not weakly irreducible and the
nonnegative eigenpairs of $T$ are
\begin{eqnarray*}
(4,(0,0,1)^T)\quad \mbox{\rm and}\quad (4.8730,(0.4365,0.5635,0)^T).
\end{eqnarray*}
Hence, we get the spectral radius of $T$ is
$4.8730$ which agrees the numerical results in Table 2.
\end{Example}

From Tables 1 and 2, we have some preliminary observations:
\begin{itemize}
\item Algorithm \ref{algo3} can find the spectral radius of a general nonnegative tensor efficiently.
\item As expected, the more dense of nonzero elements of the underlying tensor,
the higher the probability of it being weakly irreducible. We note
that function {\em sprand} which is used in our test in MATLAB does
not return a full-dense matrix even with the density parameter being
$1$, so the percentage {\bf Per} for $n=3$ is not close enough to
$100$ even when ${\bf Den}=0.9$.
\item From Definition \ref{def-new}, we note that elements of $G(T)$
are the summation of so many elements of $T$. So, it is possible that
$G(T)$ is an irreducible matrix even when $T$ is a very sparse
tensor. This can be noticed from the last several rows in Table 1.
\end{itemize}

\section{Conclusions and remarks}
\setcounter{Assumption}{0}
\setcounter{Theorem}{0} \setcounter{Proposition}{0}
\setcounter{Corollary}{0} \setcounter{Lemma}{0}
\setcounter{Definition}{0} \setcounter{Remark}{0}
\setcounter{Algorithm}{0}  \setcounter{Example}{0}
\hspace{4mm} In this paper, we proposed a new class of nonnegative
tensors --- strictly nonnegative tensors, and proved that the
spectral radii for strictly nonnegative tensors are always positive.
We discussed the relationships among the seven well-conditional
classes of nonnegative tensors, and showed that the class of
strictly nonnegative tensors strictly contains the other six classes
of nonnegative tensors. We showed that a modification of the power
method in \cite{fgh} for finding the spectral radius of a
nonnegative tensor is globally R-linearly convergent for weakly
irreducible nonnegative tensors. Based on the convergent power
method for weakly irreducible nonnegative tensors, we proposed an
algorithm for finding the spectral radius of a general nonnegative
tensor. The preliminary numerical results addressed the
effectiveness of the method.

Some issues deserve to be further investigated, such as, how to
investigate the related topics for some classes of tensors {\it
beyond nonnegative tensors}, and how to find the smallest eigenvalue
of these classes of tensors?

{\bf Acknowledgement.} We are grateful to Professor Qin Ni for inspiring Theorem \ref{add-thm-1}.

 \begin{table}[ht]\label{tab}
  \caption{Numerical results}
  \begin{center}
  \tabcolsep 4.8pt
  \renewcommand\arraystretch{1.2}
      \begin{tabular}[c]{  c || c || c c c c c || c}
      \hline \hline
                          n  & Den & $\rho$ & Per&Ite&Blks&Res&TolCpu\\ \hline \hline

                          3&      0.10&      0.148&      0.00    & 6.88  & 2.88 & 2.0516e-008&      0.30 \\\hline
                          3&      0.20&      0.409&     14.00    & 17.88  & 2.52 & 1.3521e-007&      0.38\\\hline
                          3&      0.30&      0.671&     36.00    & 25.32  & 2.06 & 6.2038e-008&      0.44\\\hline
                          3&      0.40&      0.787&     44.00    & 27.6  & 1.9 & 1.4773e-006&      0.41\\\hline
                          3&      0.50&      1.113&     66.00    & 31.06  & 1.38 & 8.1021e-008&      0.42\\\hline
                          3&      0.60&      1.032&     72.00    & 33.76  & 1.44 & 9.8837e-008&      0.44\\\hline
                          3&      0.70&      1.196&     78.00    & 31.4  & 1.34 & 8.6830e-008&      0.47\\\hline
                          3&      0.80&      1.244&     80.00    & 37.82  & 1.28 & 7.3281e-007&      0.50\\\hline
                          3&      0.90&      1.528&     86.00    & 33.26  & 1.22 & 2.1994e-006&      0.44\\\hline

                          4&      0.10&      0.368&     14.00    & 19.56  & 3.26 & 3.5644e-008&      0.47\\\hline
                          4&      0.20&      0.554&     24.00    & 27.04  & 2.74 & 7.1762e-008&      0.48\\\hline
                          4&      0.40&      1.060&     62.00    & 42.28  & 1.64 & 1.6366e-006&      0.50\\\hline
                          4&      0.80&      1.945&     90.00    & 35.76  & 1.1 & 6.5548e-008&      0.47\\\hline

                          10&      0.05&      0.466&     14.00    & 51.6  & 5.72 & 1.8168e-006&      0.94\\\hline
                          10&      0.10&      1.214&     56.00    & 49.44  & 2.34 & 5.5063e-007&      1.05\\\hline
                          10&      0.15&      2.363&     78.00    & 40.96  & 1.26 & 1.5255e-008&      0.98\\\hline
                          10&      0.20&      3.124&     88.00    & 31.44  & 1.2 & 1.2771e-008&      0.95\\\hline

                          20&      0.05&      2.537&     56.00    & 39.54  & 2.08 & 4.8818e-009&      3.58\\\hline
                          20&      0.10&      5.606&     86.00    & 31.34  & 1.14 & 4.7892e-009&      3.28\\\hline

                          30&      0.05&      6.103&     80.00    & 31.84  & 1.28 & 2.9270e-009&      9.42\\\hline
                          30&      0.10&     12.173&     84.00    & 27.06  & 1.16 & 2.6974e-009&      8.02\\\hline

                          40&      0.05&     10.698&     82.00    & 28.18  & 1.26 & 1.8279e-009&     18.02\\\hline
                          50&      0.05&     16.740&     94.00    & 27.2  & 1.06 & 1.3807e-009&     36.80\\\hline
      \hline

       \end{tabular}
  \end{center}
 \end{table}

 \begin{table}[ht]\label{tab-1}
  \caption{Numerical results for Example \ref{exm-6}}
  \begin{center}
  \tabcolsep 4.8pt
  \renewcommand\arraystretch{1.2}
      \begin{tabular}[c]{ c|| c || c || c || c || c}
      \hline \hline
                          Blk& Ite&$\alpha(x^{(k)})$&$\beta(x^{(k)})$&$\alpha(x^{(k)})-\beta(x^{(k)})$&$\|T(x^{(k)})^2-\frac{\alpha(x^{(k)})+\beta(x^{(k)})}{2}(x^{(k)})^{[2]}\|$\\ \hline \hline
1&1&     6.000&      4.000& 2.000e+000& 1.414e+000\\
 1& 2&     5.200&      4.571& 6.286e-001& 1.134e-001\\
  1&3&     4.974&      4.774& 2.002e-001& 3.573e-002\\
  1&4&     4.905&      4.841& 6.383e-002& 1.143e-002\\
  1&5&     4.883&      4.863& 2.035e-002& 3.641e-003\\
  1&6&     4.876&      4.870& 6.491e-003& 1.162e-003\\
  1&7&     4.874&      4.872& 2.070e-003& 3.704e-004\\
  1&8&     4.873&      4.873& 6.602e-004& 1.181e-004\\
  1&9&     4.873&      4.873& 2.106e-004& 3.767e-005\\
  1&10&     4.873&      4.873& 6.715e-005& 1.201e-005\\
  1&11&     4.873&      4.873& 2.141e-005& 3.832e-006\\
  1&12&     4.873&      4.873& 6.830e-006& 1.222e-006\\
  1&13&     4.873&      4.873& 2.178e-006& 3.897e-007\\
  1&14&     4.873&      4.873& 6.946e-007& 1.243e-007\\\hline\hline
  2&1&     4.000&      4.000& 0.000e+000& 0.000e+000\\\hline
   \hline

       \end{tabular}
  \end{center}
 \end{table}

\end{document}